\documentclass[12pt,thmsa]{article}
\usepackage{amsfonts}
\usepackage{amssymb}
\newtheorem{teo}{Theorem}[section]

\newtheorem{cor}[teo]{Corollary}
\newtheorem{prop}[teo]{Proposition}

\newtheorem{obs2}[teo]{Remark}

\newtheorem{tea}{Theorem}[subsection]

\newtheorem{no2}[teo]{Note}

\newtheorem{no3}[tea]{Note}

\newcommand{\Q}{{\mathbb{Q} }}

\newcommand{\GL}{{\rm GL}}

\title{Existence of families of Galois
representations and new cases of the Fontaine-Mazur conjecture}
\author{Luis V. Dieulefait\\ Centre de Recerca Matem{\'a}tica\\ Bellaterra,
Barcelona
  \thanks{
e-mail: LDieulefait@crm.es.
Supported by MECD postdoctoral fellowship at the Centre de Recerca Matem{\'a}tica from
the Ministerio de Educaci{\'o}n y Cultura
}}

\begin{document}

\maketitle
\begin{abstract}
 In a previous article, we have proved a result asserting the existence
 of a compatible family of Galois representations containing a given
 crystalline irreducible odd two-dimensional representation. We apply this
 result to establish new cases of the Fontaine-Mazur conjecture, namely,
 an irreducible Barsotti-Tate $\lambda$-adic $2$-dimensional Galois
 representation unramified at $3$ and such that the traces $a_p$
 of the images of Frobenii verify $\Q(\{ a_p^2 \}) = \Q $ always
 comes from an abelian variety. We also show the non-existence of
 irreducible Barsotti-Tate $2$-dimensional Galois representations of conductor $1$
and apply this to the irreducibility of Galois representations on
level 1 genus 2 Siegel cusp forms.

\end{abstract}

\section{Existence of families}
The following is a slight generalization of a result proved in [D3],
which follows from the results and techniques in [T1], [T2]:
 \begin{teo}
  \label{teo:Taylor} Let q be an odd prime, and $Q$ a place above $q$. Let
  $\sigma_{Q}$ be a two dimensional odd irreducible $Q$-adic Galois
  representation (of the absolute Galois group of $\Q$, continuous)
  ramified only at $q$ and at a finite set of primes $S$. Assume that
  $\sigma_Q$ is crystalline at $q$, with Hodge-Tate weights  $\{0, w \}$ ($w$ odd).
  Assume also that $q \geq 2w + 1$.
  Then, there exists a  compatible family of
  Galois representations $\{ \sigma_{\lambda} \}$ containing $
  \sigma_{Q}$, such that
   for every $\ell \not\in
  S$, $\lambda \mid \ell$, the representation $\sigma_{\lambda}$ is
  unramified outside $\{ \ell \} \cup S$ and is crystalline at $\ell$ with
  Hodge-Tate weights $\{ 0 , w \}$. Moreover, the family $\{ \sigma_{\lambda} \}$
   is strictly compatible (see [T1] for the definition) and all its members are
    irreducible.
  \end{teo}
  Remark: the proof given in [D3] applies in this generality, just
  observe that for the case of $w =1$, in general when Taylor proves
  that the restriction of $\sigma_Q$ to some totally real field $F$ will correspond
  to a representation attached to a Hilbert modular form $h$, this is not
  enough to conclude that the representation is motivic (the construction
   of Blasius-Rogawski does not apply in some cases), thus is not
   enough in general to prove the Fontaine-Mazur conjecture, but in any
   case there is a strictly compatible family of Galois representation
   attached to $h$ (by previous results of Taylor) and so the argument
   in [D3] (descent of this family to a compatible family of $G_\Q$-representations)
   can be applied also in this case. The family obtained will be
   strictly compatible, as follows from the results in [T1]. All representations
    in the family are irreducible because when restricted to the Galois group of
    $F$ they agree with the modular Galois representations attached to
     $h$, which are irreducible (because $h$ is cuspidal).\\

\section{Fontaine-Mazur for ``projectively rational" Barsotti-Tate representations}
Now suppose that we are given a representation $\sigma_Q$ as in
theorem \ref{teo:Taylor} with $w = 1$ and $\det (\sigma_Q) = \chi$.
Then, we will prove the following:
\begin{teo}
\label{teo:abel} Assume that $\sigma_Q$ verifies also the following two
conditions:\\
1) If $q \neq 3$, then $3 \not\in S$ ($\sigma_Q$ unramified at $3$).\\
2) The traces $\{ a_p \}$ of the images of Frobenii,
 for every $p \neq q, p \not\in S$
verify: $a_p^2 \in \mathbb{Z}$.\\
Then $\sigma_Q$ can be attached to an
abelian variety $A$ (of course, $A$ is of $\GL_2$-type). \\
\end{teo}
Remark: In particular, in case all the traces are integers, this result
shows that a Galois representation that ``looks like" the one attached
to an elliptic curve with good reduction at $3$
does indeed come from such an elliptic curve.\\

Remark: A similar result is proved in [T1] without restriction on the traces
  but (sticking to the case $w=1$) with the extra assumption
  that  there is
 a prime $u \in
S$ such that the restriction of $\sigma_Q$ to the decomposition group $D_u$
 is of a particular type
 (corresponding to discrete
 series under the local Langlands correspondence). Thus, this result
  does not apply, for example, if $\sigma_Q$ has conductor 1 or
 is semistable (i.e., unipotent) locally at every prime of $ S$. Therefore, in the
  semistable case, the results of [T1] are not enough to prove the Fontaine-Mazur
   conjecture for $\sigma_Q$.\\

Proof:
   The proof follows from the combination of theorem
   \ref{teo:Taylor} with modularity results … la Wiles. We know that
    there exists a strictly compatible family $\{ \sigma_\lambda \}$
   containing $\sigma_Q$. Take $t \mid 3$ and consider the Galois
   representation $\sigma_t$ (if $q=3$, just take $t=Q$).
   As $3 \not\in S$, we know that $\sigma_t$
    is crystalline at $t$, and has Hodge-Tate weights $\{ 0, 1\}$.
    Following the initial idea of Wiles (see also [D1]) we know, from
    condition 2) in the theorem (via results of Langlands and Tunnell), that the
    residual representation $\bar{\sigma}_t$ will be either modular or
    reducible. The information we have on $\sigma_t$ is enough then to
    conclude, via a combination of modularity results … la Taylor-Wiles
    and Skinner-Wiles, that $\sigma_t$ is modular. This is a non-trivial
    assertion, but this is done in [D1]
    in exactly the same situation! So we conclude that the family $\{\sigma_\lambda\}$
    is modular, and from $w=1$ we easily check that it will be attached
    to a weight $2$ cusp form f. This proves that $\sigma_Q$ can be
    attached to the abelian variety $A_f$. \\

    Remark: It follows from condition 2) in the theorem that the
    variety $A_f$ will have a large endomorphism algebra (cf. [R]).\\

\section{Non-existence of Barsotti-Tate representations of conductor 1}
In this section we will prove the following result:
\begin{teo}
\label{teo:nohay} There are no
irreducible Barsotti-Tate $2$-dimensional Galois representations of conductor $1$
and odd residual characteristic.
 \end{teo}
Remark: By Barsotti-Tate we mean crystalline Galois representations as
$\sigma_Q$ in theorem \ref{teo:Taylor} with $w=1$.\\

Proof: As in the previous section, $q$ is odd and $\sigma_Q$ has $w=1$.
Now there is no restriction in the field of coefficients, but $S$ is
empty. If $q \neq 3$, once again we use the results in section 1 to
construct a strictly compatible family $\{\sigma_\lambda\}$ containing
$\sigma_Q$. We consider again $\sigma_t$ for $t \mid 3$ (just take $t=Q$ if
$q=3$). From strict compatibility, this $t$-adic Galois representation
will be unramified outside 3. But such a representation can not exist, as
was proved in [D2]. Let us briefly recall the argument for the convenience
of the reader: the residual representation $\bar{\sigma}_t$
 has coefficient in a finite
field of characteristic $3$ and is unramified outside $3$. Then, a
result of Serre tells us that it must be reducible. An application of
results of Skinner-Wiles (that can be applied because $\sigma_t$ is Barsotti-Tate)
shows that $\sigma_t$ is modular, but this is a contradiction because
it must correspond to a level 1, weight 2, cuspidal modular form.\\

\begin{cor}
\label{teo:masqueirred} Let $f$ be a genus $2$, level $1$, cuspidal
Siegel modular form (Hecke eigenform) having multiplicity one and weight $k>3$. Suppose
that $f$ is not a Maass spezialform. Then, for every odd prime $\ell$,
$\lambda \mid
\ell$ in $E$= field generated by the eigenvalues of $f$, the Galois
representation $\rho_{f,\lambda}$ attached to $f$ is absolutely
irreducible. In particular, this representation can be defined over
$E_\lambda$.
\end{cor}
Proof: In [D2] it is shown that the only possible reducible case is:
$\rho_{f,\lambda} \cong \sigma_{1,\lambda} \oplus \sigma_{2,\lambda}$
where one of the two (necessarily irreducible) components, say
$\sigma_{2,\lambda}$, is crystalline, with Hodge-Tate weights
 $\{k-2 , k-1\}$, and has conductor 1. But theorem \ref{teo:nohay} says
  that such a representation (after twisting by $\chi^{2-k}$)
   can not exist. This proves the
  corollary.\\

  Remark: This improves the main result of [D3], for the level
  1 case: ``uniformity of reducibility" now does not make sense, because all
  representations will be irreducible.\\

   In the semistable case, the
  main theorem of [D3] can be extended, just by applying theorem
  \ref{teo:Taylor}, to be valid for every prime:

  \begin{prop}
  \label{teo:lomas} Let $\{ \rho_\lambda\}$ be a compatible family of
  4-dimensional, pure, symplectic, Galois representations, with
  finite ramification set $S$, semistable (at every place in $S$) in the sense of [D3]
   and such that for every $\ell \not\in S$, $\lambda \mid \ell$,
  $\rho_\lambda$ is crystalline with Hodge-Tate weights $\{0 , k-2, k-1,
  2k-3\}$. Then, if for some $q >2, q\not\in S$, $Q \mid q$, the representation
  $\rho_Q$ is reducible,
    all the representations in the family are
  reducible.
  \end{prop}

Proof: It is known (cf. [D2], [D3]) that the only possible reducible case
 is: $\rho_Q \cong \sigma_{1,Q} \oplus \sigma_{2,Q}$, where the two components
  are irreducible odd two-dimensional Galois representations of the same
  determinant, one of them having Hodge-Tate weights $\{0, 2k-3\}$, the
  other having weights $\{ k-2 , k-1\}$.
   If $q \geq 4k-5$, an application of theorem \ref{teo:Taylor} to both
   components (for one of them, you may twist by a power of $\chi$ before applying
    the theorem, and then untwist to obtain the desired family)
     proves the result, because from compatibility and
   Cebotarev density theorem it is clear that the families $\{\sigma_{1,\lambda}\}$
    and $\{\sigma_{2,\lambda} \}$ verify: $\rho_{\lambda} \cong
     \sigma_{1,\lambda} \oplus \sigma_{2,\lambda}$ (up to
     semisemplification) for every $\lambda$.\\
     If $2 <q < 4k-5$, then we can only apply (with the twist and untwist trick)
      theorem \ref{teo:Taylor} to one of the components, say $\sigma_{2,
      Q} \in \{ \sigma_{2,\lambda }\}$. The techniques of [D3]
      apply precisely in this situation: it is enough (because we are
      assuming semistability) to have the ``existence of a family" result
       for one component, to conclude reducibility of $\{\rho_\lambda\}$
       for almost every prime
      But this implies that we can take now a
       second prime $r$ as large as we want ($r \geq 4k-5, r \not\in S$), and $R \mid
       r$ such that $\rho_R$ is reducible: $\rho_R \cong \sigma_{1,R} \oplus
        \sigma_{2,R}$, and now because $r$ is sufficiently large we can
        apply theorem \ref{teo:Taylor} also to the other component,
        and conclude  as before that the whole family $\{ \rho_\lambda\}$
        is reducible.\\

Remark: When we assume that $\rho_Q$ is reducible, reducibility must
necessarily ``occur over the field of coefficients" (the field generated by the traces
 of the images of Frobenii), i.e., the field of coefficients of
 $\rho_Q$ must contain the fields of coefficients of its irreducible
 components (this is required in the proof above, because we have used
  the results of [D3]). That reducibility always occurs over the field of coefficients
   was proved in [D2] for $q \geq 4k-5$ but the proof
 holds for every odd prime $q$: if we assume that this is not the case, following
  the arguments in [D2], we conclude that the field of coefficients
 of $\sigma_{2,Q}$ is an infinite extension of that of $\rho_Q$ (to orient the
 reader, let us point out that this result
 strongly depends on the fact that $\rho_Q$ has four different Hodge-Tate
 weights, and that when reducibility is not over the fields of coefficients this
 forces the
 images of the representations $\{ \rho_\lambda \}$ to be ``generically
 large", cf. [D2], [D3])), but
 this contradicts the fact that $\sigma_{2,Q}$ (after twisting) is
 potentially modular (thus, the field of coefficients of $\rho_Q |_{G_F}$ is a number
  field for some totally real number field $F$,
   and a fortiori the field of coefficients of $\rho_Q$ is also a finite extension
    of $\Q$).

\section{Bibliography}



[D1] Dieulefait, L., {\it Modularity of Abelian Surfaces with Quaternionic
Multiplication}, Math. Res. Letters {\bf 10} no. 2-3 (2003)\\
\newline
[D2] Dieulefait, L.,  {\it Irreducibility of Galois actions on level $1$
Siegel cusp forms}, preprint, (2002)\\
available at http://www.math.leidenuniv.nl/gtem/view.php (preprint
19)\\
\newline
[D3] Dieulefait, L., {\it Uniform behavior  of families of Galois
representations on Siegel modular forms}, preprint, (2002)\\
available at http://www.math.leidenuniv.nl/gtem/view.php (preprint
68)\\
\newline
[R] Ribet, K., {\it Abelian varieties over $\Q$ and modular forms},
Algebra and Topology 1992, KAIST Math. Workshop, Taejon, Korea
(1992) 53-79\\
\newline
[T1] Taylor, R., {\it On the meromorphic continuation of degree two
 L-functions}, preprint, (2001);
 available at http://abel.math.harvard.edu/$\sim$rtaylor/ \\
\newline
[T2] Taylor, R., {\it Galois Representations}, proceedings of ICM, Beijing,
August 2002, longer version available at
http://abel.math.harvard.edu/$\sim$rtaylor/ \\
\newline



\end{document}